\documentclass[12pt,reqno]{amsproc}
\usepackage{amsthm,amsfonts,amssymb,amscd,euscript}
\newcommand{\Diff}[2]{\operatorname{Diff}_{#1}{\left(#2\right)}}
\newcommand{\Dif}[1]{\operatorname{Diff}_{#1}}

\newcommand{\NE}{\overline{NE}}
\newcommand{\Supp}{\operatorname{Supp}}
\newcommand{\id}{\operatorname{id}}
\newcommand{\Exc}{\operatorname{Exc}}

\newcommand{\red}{\operatorname{red}}
\newcommand{\comment}[1]{}
\newcommand{\ep}{\varepsilon}
\renewcommand{\bar}[1]{\overline{#1}}
\newcommand{\Sing}[1]{\operatorname{Sing}( #1)}
\newcommand{\QQ}{{\mathbb Q}}
\newcommand{\ZZ}{{\mathbb Z}}
\newcommand{\NN}{{\mathbb N}}
\newcommand{\CC}{{\mathbb C}}
\newcommand{\PP}{{\mathbb P}}
\newcommand{\FF}{{\mathbb F}}
\newcommand{\OOO}{{\mathcal O}}
\newcommand{\KKK}{{\EuScript{K}}}
\newcommand{\down}[1]{\left\lfloor #1\right\rfloor}

\renewcommand{\emptyset}{\varnothing}
\newtheorem{theorem}[subsection]{Theorem}
\newtheorem*{theorem*}{Theorem}
\newtheorem{lemma}[subsection]{Lemma}
\newtheorem{claim}[subsection]{Claim}
\newtheorem{proposition}[subsection]{Proposition}

\newtheorem{conjecture}[subsection]{Conjecture}
\theoremstyle{definition}
\newtheorem{definition}[subsection]{Definition}
\newtheorem{example}[subsection]{Example}
\title{On classification of Mori contractions: Elliptic curve case}
\author{Yuri~G.~Prokhorov}
\thanks{This work was partially supported by grants
INTAS-OPEN-97-2072 and RFFI 99-01-01132}
\email{prokhoro@mech.math.msu.su}
\address{Department of Mathematics (Algebra Section),
Moscow University, 117234 Moscow, Russia}

\begin{document}
\begin{abstract}
We study Mori's three-dimensional contractions $f\colon X\to Z$.
It is proved that on the ``good'' model $(\bar X,\bar S)$ there
are no elliptic components of $\Dif{\bar S}$ with coefficients
$\ge 6/7$.
\end{abstract}
\maketitle

\section{Introduction}
\subsection{}
\label{notation}
Let $f\colon X\to Z\ni o$ be an extremal log terminal contraction
over $\CC$, that is:
\par\medskip\noindent
$X$ is a normal algebraic $\QQ$-factorial threefold with at worst
log terminal singularities, $f$ is a projective morphism such that
$f_*\OOO_X=\OOO_Z$, $\rho(X/Z)=1$ and $-K_X$ is $f$-ample.
\par\medskip\noindent
We assume that $\dim(Z)\ge 1$ and regard $(Z\ni o)$ as a
sufficiently small Zariski neighborhood. Such contractions
naturally appear in the Minimal Model Program \cite{KMM}. By
$\Exc(f)\subset X$ denote the exceptional locus of $f$.

According to the general principle introduced by Shokurov
\cite{Sh} all such contractions can be divided into two classes:
\emph{exceptional} and \emph{nonexceptional}. A contraction
$f\colon X\to Z\ni o$ such as in \ref{notation} is said to be
exceptional if for any complement $K_X+D$ near $f^{-1}(o)$ there
is at most one divisor $E$ of the function field $\KKK(X)$ with
discrepancy $a(E,D)=-1$. The following is a particular case of the
theorem proved in \cite{Sh1} and \cite{P} (see also \cite{PSh}).
\begin{theorem}
Notation as above. Assume that $f\colon X\to Z\ni o$ is
nonexceptional. Then for some $n\in\{1,2,3,4,6\}$ there is a
member $F\in |-nK_X|$ such that the pair $\left(X,\frac1nF\right)$
is log canonical near $f^{-1}(o)$.
\end{theorem}
Thus nonexceptional contractions have a ``good'' member in
$|-nK_X|$, $n\in\{1,2,3,4,6\}$. The most important case is the
case of \emph{Mori contractions}, i.~e.  when $X$ has only
terminal singularities:

\begin{conjecture}
Notation as in \ref{notation}. Assume that $X$ has at worst
terminal singularities. Then $f\colon X\to Z\ni o$ is
nonexceptional.
\end{conjecture}

Similar to the classification of three-dimensional terminal
singularities, this fact should be the key point in the
classification Mori contractions. For example it is very helpful
in the study of three-dimensional flips \cite{K}, \cite{Mo},
\cite{Sh}.

Methods of \cite{Sh1}, \cite{P}, \cite{PSh} use inductive
procedure of constructing divisors in $|-nK|$. This procedure
works on so-called \textit{good} model of $X$ over $Z$. Roughly
speaking, a good model is a birational model $\bar Y$ equipped
with a prime divisor $\bar S$ such that the pair $(\bar Y,\bar S)$
is plt and $-(K_{\bar Y}+\bar S)$ is nef and big over $Z$.

If $f$ is exceptional, then $\bar S$ is a projective surface.
Adjunction Formula \ref{def-adj} gives us that $(\bar S,\Dif{\bar
S})$ is a \textit{klt log del Pezzo surface}. Moreover,
exceptionality of $f$ implies that the projective log pair $(\bar
S,\Dif{\bar S})$ is \textit{exceptional}, by definition this means
that any complement $K_{\bar S}+\Dif{\bar S}^+$ is klt
\cite[Prop.~ 2.4]{P}. Thus our construction gives the following
correspondence:

\[
\begin{array}{ccc}
\text{\begin{tabular}{|l|}\hline Exceptional contractions \\
$f\colon X\to Z$ as in \ref{notation}\\ \hline
\end{tabular}}&\longrightarrow&
\text{\begin{tabular}{|l|}\hline Exceptional log del Pezzos \\
$(\bar S,\Delta=\Dif{\bar S})$ \\ \hline\end{tabular}}\\
\end{array}
\]

\subsection{}
For exceptional log del Pezzos $(\bar S,\Delta)$ Shokurov
introduced the following invariant:
\[
\begin{array}{ll}
\delta=\delta(\bar S,\Delta)=&\text{number of divisors ${E}$ of
$\KKK(\bar S)$}\\ &\text{with discrepancy $a(E,\Delta)\le
-6/7$.}\\
\end{array}
\]
He proved that $\delta\le 2$, classified log surfaces with
$\delta=2$ and showed that in the case $\delta=1$ the (unique)
divisor $E$ with $a(E,\Delta)\le -6/7$ is represented by a curve
of arithmetical genus $\le 1$ (see \cite{Sh1}, \cite{Pr-lect}).

The aim of this short note is to exclude the case of Mori
contractions with $\delta=1$ and elliptic curve $E$:
\begin{theorem}
\label{main}
Notation as in \ref{notation}. Assume that ${\delta}(\bar
S,\Dif{\bar S})=1$. Write $\Dif{\bar S}=\sum \delta_i\bar
\Delta_i$, where $\bar \Delta_i$ are irreducible curves. If
$\delta_{i_0}\ge 6/7$ for some $i_0$, then $p_a(\bar
\Delta_{i_0})=0$.
\end{theorem}
The following example shows that Theorem~\ref{main} cannot be
generalized to the klt case.
\begin{example}[\cite{IP}]
Let $(Z\ni o)$ be the hypersurface canonical singularity
$x_1^2+x_3^3+x_3^{11}+x_4^{12}=0$ and let $f\colon X\to Z$ be the
weighted blowup with weights $(66,44,12,11)$. Then $f$ satisfies
conditions of \ref{notation} and we have case~\ref{plt}. It was
computed in \cite{IP} that $S$ is the weighted projective plane
$\PP(3,2,1)$ and $\Dif{S}=\frac{10}{11}C+\frac{1}{2}L$, where $C$
is an elliptic curve.
\end{example}

Log del Pezzo surfaces of elliptic type (like $(\bar S,\Dif{\bar
S})$ in the above theorem) were classified by T.~Abe \cite{Abe}.
Our proof uses different, very easy arguments.

\section{Preliminary results}
In this paper we use terminologies of Minimal Model Program
\cite{KMM}, \cite{Ut}. For the definition of \textit{complements}
and their properties we refer to \cite[Sect. 5]{Sh}, \cite[Ch.
19]{Ut} and \cite{Pr-lect}.

\begin{definition}[{\cite[Sect. 3]{Sh}}, {\cite[Ch. 16]{Ut}}]
\label{def-adj}
Let $X$ be a normal variety and let $S\ne\emptyset$ be an
effective reduced divisor on $X$. Let $B$ be a $\QQ$-divisor on
$X$ such that $S$ and $B$ have no common components. Assume that
$K_X+S$ is lc in codimension two. Then the \textit{different} of
$B$ on $S$ is defined by
\begin{equation*}
K_S+\Diff{S}{B}\equiv (K_X+S+B)|_S.
\end{equation*}
Usually we will write simply $\Dif{S}$ instead of $\Diff{S}0$.
\end{definition}

\begin{theorem}[Inversion of Adjunction {\cite[17.6]{Ut}}]
\label{Inv-Adj}
Let $X$ be a normal variety and let $D$ be a boundary on $X$.
Write $D=S+B$, where $S=\down{D}$. Assume that $K_X+S+B$ is
$\QQ$-Cartier. Then $(X,S+B)$ is plt near $S$ iff $S$ is normal
and $(S,\Diff{S}{B})$ is klt.
\end{theorem}

\begin{definition}[\cite{Pr1}]
\label{def-plt-b}
Let $X$ be a normal variety and let $g\colon Y\to X$ be a
birational contraction such that the exceptional locus of $g$
contains exactly one irreducible divisor, say $S$. Assume that
$K_Y+S$ is plt and $-(K_Y+S)$ is $f$-ample. Then $g\colon
(Y\supset S)\to X$ is called a \textit{plt blowup} of $X$.
\end{definition}

The key point in the proof of Theorem~\ref{main} is the following
proposition.
\begin{proposition}
\label{prop-main}
Let $(X\ni P)$ be a three-dimensional terminal singularity and let
$g\colon (Y,S)\to X$ be a plt blowup with $f(S)=P$. Write
$\Dif{S}=\sum \delta_i\Delta_i$, where $\Delta_i$ are irreducible
curves, and assume that $\delta_0\ge 6/7$ for some $i_0$. Further,
assume that $S$ is smooth at singular points of $\Delta_{i_0}$.
Then $p_a(\Delta_{i_0})=0$.
\end{proposition}

\begin{lemma}[cf. {\cite[Cor. 5]{Pr1}}]\label{first}
Let $(X\ni P)$ be a three-dimensional terminal singularity and let
$g\colon (Y,S)\to X$ be a plt blowup with $f(S)=P$. Then there is
a boundary $\Upsilon\ge \Dif{S}$ on $S$ such that
\begin{enumerate}
\item
$\down \Upsilon\ne 0$;
\item
$-(K_S+\Upsilon)$ is ample.
\end{enumerate}
Moreover, $K_S+\Dif{S}$ has a non-klt $1$, $2$, $3$, $4$, or
$6$-complement.
\end{lemma}
\begin{proof}
Regard $(X\ni P)$ as an analytic germ. It was shown in
\cite[Sect.~ 6.4]{YPG} that the general element $F\in |-K_X|$ has
a normal Du Val singularity at $P$. By Inversion of Adjunction
\ref{Inv-Adj}, $K_X+F$ is plt. Consider the crepant pull-back
$K_Y+aS+F_Y= f^*(K_X+F)$, where $F_Y$ is the proper transform of
$F$ and $a<1$. Since both $K_Y+S$ and $g^*K_X$ are $\QQ$-Cartier,
so are $S$ and $F_Y$. Clearly, $-(K_Y+S+F_Y)$ is $f$-ample. Let
$\Upsilon':= \Diff{S}{F_Y}$. Then $\down \Upsilon'\ne 0$ and
$-(K_S+\Upsilon')$ is ample. Therefore $\Upsilon$ can be found in
the form $\Upsilon=\Diff{S}{F_Y}+t(\Upsilon'-\Diff{S}{F_Y})$ for
suitable $0<t\le 1$.

Take $\Delta=\Diff{S}{F_Y}+\lambda(\Upsilon-\Diff{S}{F_Y})$ for
$0<\lambda\le 1$ so that $K_S+\Delta$ is lc but not klt (and
$-(K_S+\Delta)$ is ample). By \cite[Sect.~2]{Sh1} (see also
\cite[5.4.1]{Pr-lect}) there exists either an $1$, $2$, $3$, $4$,
or $6$-complement of $K_S+\Delta$ which is not klt.
\end{proof}

A very important problem is to prove the last lemma without using
\cite[Sect. 6.4]{YPG}, i.e. the classification of terminal
singularities. This can be a way in higher-dimensional
generalizations.

\begin{proof}[Proof of Proposition~{\rm \ref{prop-main}}]
Put $C:=\Delta_{i_0}$ and let $\delta_0=1-1/m$, $m\ge 7$. Assume
that $p_a(C)\ge 1$. Let $\Upsilon$ be such as in Lemma~\ref{first}
and let $K_S+\Theta$ be a non-klt $1$, $2$, $3$, $4$, or
$6$-complement of $K_S+\Dif{S}$. Using that the coefficients of
$\Dif{S}$ are standard \cite[Prop. 3.9]{Sh} it is easy to see that
$\Theta\ge \Dif{S}$ and $\Theta\ge C$ \cite[Sect. 4.7]{Pr-lect}.
In particular, $K_S+C$ is lc.

Further, $C\not \subset\down{\Upsilon}$. Indeed, otherwise by
Adjunction we have
\[
-\deg K_C\ge -\deg(K_C+\Diff{C}{\Upsilon-C})=-(K_S+\Upsilon)\cdot
C>0
\]
This implies $p_a(C)=0$, a contradiction.

By Lemma~\ref{wheel-1} below, $\Theta=C$, $p_a(C)=1$, $S$ is
smooth along $C$ and has only Du Val singularities outside.
Therefore, $\Dif{S}=(1-1/m)C$ and $-K_S\equiv C\equiv
-m(K_S+(1-1/m)C)$ is ample (see \ref{def-plt-b}). Thus $S$ is a
del Pezzo surface with at worst Du Val singularities. Since
$C\not\subset\down\Upsilon$, we can write $\Upsilon=\alpha
C+L+\Upsilon^{o}$, where $L$ is an irreducible curve, $1>\alpha\ge
1-1/m\ge 6/7$, $C\not \subset\Supp(\Upsilon^{o})$ and
$\Upsilon^{o}\ge 0$. Further,
\begin{multline}
\label{eq-Up}
0<K_S\cdot (K_S+\Upsilon)=K_S\cdot (K_S+\alpha C+L+\Upsilon^{o})\\
\le K_S\cdot ((1-\alpha)K_S+L) \le \frac17K_S^2+K_S\cdot L.
\end{multline}
Thus $K_S^2>-7K_S\cdot L\ge 7$.

Let $S^{\min}\to S$ be the minimal resolution. By Noether's
formula, $K_{S}^2+\rho(S^{\min})= K_{S^{\min}}^2+\rho(S^{\min})=
10$. Thus, $8\le K_S^2\le 9$ and $\rho(S^{\min})\le 2$. In
particular, $S$ either is smooth or has exactly one singular point
which is of type $A_1$. By \eqref{eq-Up}, $-K_S\cdot L=1$. Similar
to \eqref{eq-Up} we have
\[
0<-(K_S+\Upsilon)\cdot L\le -(1-\alpha)K_S\cdot L -L^2.
\]
Hence $L^2<1-\alpha\le 1/7$, so $L^2\le 0$. This means that the
curve $L$ generates an extremal ray on $S$ and $\rho(S)=2$.
Therefore, $S$ is smooth and $K_S^2=8$. In this case, $S$ is a
rational ruled surface ($\PP^1\times\PP^1$ or $\FF_1$). Let $\ell$
be a general fiber of the rulling. Then
\[
0<-(K_S+\alpha C+L+\Upsilon^{o})\cdot \ell\le -(1-\alpha)K_S\cdot
\ell-L\cdot \ell\le \frac27-L\cdot \ell.
\]
Hence $L\cdot \ell=0$, so $L$ is a fiber of $S\to \PP^1$ and
$-K_S\cdot L=2$. This contradicts \eqref{eq-Up}.
\end{proof}

\begin{lemma}[see {\cite[Lemma 8.3.6]{Pr-lect}}]
\label{wheel-1}
Let $(S,C+\Xi)$ be a rational projective log surface, where $C$ is
the reduced and $\Xi$ is an arbitrary boundary. Assume that
$K_S+C+\Xi$ is lc, $S$ is smooth at singular points of $C$,
$K_S+C+\Xi \equiv 0$, $C$ is connected and $p_a(C)\ge 1$. Then
$\Xi =0$, $K_S+C\sim 0$, $S$ is smooth along $C$ and has only
DuVal singularities outside.
\end{lemma}
\begin{proof}
By Adjunction, $K_{C}+\Diff{C}{\Xi}=0$. Hence $\Diff{C}{\Xi}=0$.
This shows that $C\cap\Supp(\Xi)=\emptyset$ and $S$ has no
singularities at points of $C\setminus\Sing{C}$ (see \cite[Prop.
16.6, Cor. 16.7]{Ut}). Let $\mu\colon \tilde S\to S$ be the
minimal resolution and let $\tilde C$ be the proper transform of
$C$ on $\tilde S$. Define $\tilde \Xi$ as the crepant pull-back:
$K_{\tilde S}+\tilde C+\tilde \Xi=\mu^*(K_S+C+\Xi)$. It is
sufficient to show that $\tilde\Xi=0$. Assume the converse.
Replace $(S,C+\Xi)$ with $(\tilde S, \tilde C+\tilde \Xi)$. It is
easy to see that all the assumptions of the lemma holds for this
new $(S,C+\Xi)$. Contractions of $(-1)$-curves again preserve the
assumptions. Since $C$ and $\Supp(\Xi)$ are disjoint, whole $\Xi$
cannot be contracted. Thus we get $S\simeq\PP^2$ or $S\simeq\FF_n$
(a rational ruled surface). In both cases simple computations
gives us $\Xi=0$.
\end{proof}

\section{Construction of a good model}
Notation as in \ref{notation}. We recall briefly the construction
of \cite{P} (see also \cite{PSh}). Assume that $f\colon X\to Z\ni
o$ is exceptional. Let $K_X+F$ be a complement which is not klt.
There is a divisor $S$ of $\KKK(X)$ such that $a(S,F)=-1$. Since
$f$ is exceptional, this divisor is unique.

\subsection{}
\label{non-plt}
First we assume that the center of $S$ on $X$ is a curve or a
point. Then $\down F=0$. Let $g\colon Y\to X$ be a minimal log
terminal modification of $(X,F)$ \cite[17.10]{Ut}, i.~e. $g$ is a
birational projective morphism such that $Y$ is $\QQ$-factorial
and $K_Y+S+A=g^*(K_X+F)$ is dlt, where $A$ is the proper transform
of $F$. In our situation, $K_Y+S+A$ is plt. By \cite[Prop.
2.17]{Ut}, $K_Y+S+(1+\ep)A$ is also plt for sufficiently small
positive $\ep$.

\subsubsection{}\label{diagram} Run $(K+S+(1+\ep)A)$-Minimal Model
Program over $Z$:
\begin{equation}
\label{dia}
\begin{array}{ccc}
Y&\dashrightarrow &\bar Y \\ \downarrow\lefteqn{\scriptstyle{g}} &
&\downarrow\lefteqn{\scriptstyle{q}} \\X
&\stackrel{f}{\longrightarrow} &Z \\
\end{array}
\end{equation}
Note that $K_{\bar Y}+\bar S+(1+\ep)\bar A\equiv \ep \bar A\equiv
-\ep(K_{\bar Y}+\bar S)$. At the end we get so-called \emph{good
model}, i.~e. a log pair $(\bar Y,\bar S+\bar A)$ such that one of
the following holds:
\begin{enumerate}
\item[(\thesubsection.A)]
$\rho(\bar Y/Z)=2$ and $-(K_{\bar Y}+\bar S)$ is nef and big over
$Z$;
\item[(\thesubsection.B)]
$\rho(\bar Y/Z)=1$ and $-(K_{\bar Y}+\bar S)$ is ample over $Z$.
\end{enumerate}

\subsection{}
\label{plt}
Now assume that the center of $S$ on $X$ is of codimension one.
Then $S=\down F$. In this case, we put $g=\id$, $Y=X$ and $A=F-S$.
If $-(K_X+S)$ is nef, then we also put $\bar Y=X$, $\bar S=S$.
Assume $-(K_X+S)\equiv A$ is not nef. Since $A$ is effective, $f$
is birational. If $f$ is divisorial, then it must contract a
component of $\Supp(A)$. Thus $S$ is not a compact surface, a
contradiction (see \cite[Prop. 2.2]{P}). Therefore $f$ is a
flipping contraction. In this case, in diagram \eqref{dia} the map
$X=Y\dashrightarrow \bar Y$ is the corresponding flip.

Since $f$ is exceptional, in both cases \ref{non-plt} and
\ref{plt} we have $f(g(S))=q(\bar S)=o$ (see \cite[Prop.~2.2]{P}).
Adjunction Formula \ref{def-adj} gives us that $(\bar S,\Dif{\bar
S})$ is a \textit{klt log del Pezzo surface}, i.~e. $(\bar
S,\Dif{\bar S})$ is klt and $-(K_{\bar S}+\Dif{\bar S})$ is nef
and big (see \cite[Lemma~2.4]{P}). Moreover, exceptionality of $f$
implies that the pair $(\bar S,\Dif{\bar S})$ is exceptional,
i.~e. any complement $K_{\bar S}+\Dif{\bar S}^+$ is klt.

\begin{proposition}
\label{main1}
Let $f\colon X\to Z\ni o$ be as in \ref{notation}. Assume that $f$
is exceptional. Furthermore,
\begin{enumerate}
\item
if $f$ is divisorial, we assume that the point $(Z\ni o)$ is
terminal;
\item
in the case $\dim(Z)=1$, we assume that singularities of
$X\setminus f^{-1}(o)$ are canonical.
\end{enumerate}
Then case (\ref{non-plt}.B) does not occur.
\end{proposition}
\begin{proof}
Assume the converse. Then $\rho(\bar Y/Z)=1$ and $q\colon \bar
Y\to Z$ is also an exceptional contraction as in \ref{notation}.
First, we consider the case when $f$ is divisorial. Then $q$ is a
plt blowup of a terminal point $(Z\ni o)$ and $q(\bar S)=o$ (see
\cite[Prop. 2.2]{P}). By Lemma~\ref{first}, $(\bar S,\Dif{\bar
S})$ has a non-klt complement. This contradicts
\cite[Prop.~2.4]{P}. Clearly, $f$ cannot be a flipping contraction
(because, in this case, the map $Y\dashrightarrow\bar Y$ must be
an isomorphism in codimension one). If $\dim(Z)=2$, then $q$ is
not equidimensional, a contradiction.

Finally, we consider the case $\dim(Z)=1$ (and $\bar S$ is the
central fiber of $q$). Let $\bar F$ be a general fiber of $q$ (a
del Pezzo surface with at worst Du Val singularities). Consider
the exact sequence
\[
0\longrightarrow \OOO_{\bar Y}(-K_{\bar Y}-\bar F) \longrightarrow
\OOO_{\bar Y}(-K_{\bar Y})\longrightarrow \OOO_{\bar F}(-K_{\bar
F}) \longrightarrow 0
\]
By Kawamata-Viehweg Vanishing \cite[Th. 1-2-5]{KMM},
$R^1q_*\OOO_{\bar Y}(-K_{\bar Y}-\bar F)=0$. Hence there is the
surjection
\[
H^0(\bar Y,\OOO_{\bar Y}(-K_{\bar Y}))\longrightarrow H^0(\bar F,
\OOO_{\bar F}(-K_{\bar F})) \longrightarrow 0.
\]
Here $H^0(\bar F, \OOO_{\bar F}(-K_{\bar F}))\neq 0$ (because
$-K_{\bar F}$ is Cartier and ample). Therefore, $H^0(\bar
Y,\OOO_{\bar Y}(-K_{\bar Y}))\neq 0$. Let $\bar G\in |-K_{\bar
Y}|$ be any member. Take (positive) $c\in\QQ$ so that $K_{\bar
Y}+\bar S+c\bar G$ is lc and not plt. Clearly $c\le 1$, so
$-(K_{\bar Y}+\bar S+c\bar G)\equiv -(1-c)K_{\bar Y}$ is $q$-nef.
By Base Point Free Theorem \cite[Th. 3-1-1]{KMM} there is a
complement $K_{\bar Y}+\bar S+c\bar G+L$, where $nL\in |-n(K_{\bar
Y}+\bar S+c\bar G)|$ for sufficiently big and divisible $n$ and
this complement is not plt, a contradiction with exceptionality
(see \cite[Prop.~2.4]{P}).
\end{proof}

\section{Proof of Theorem~\ref{main}}
In this section we use notation and assumptions of \ref{notation}
and Theorem~\ref{main}.

If $g=\id$, then $Y=X$ and $\bar Y$ have only terminal
singularities (see \ref{plt} and \cite[5-1-11]{KMM}). Then
$\Dif{\bar S}=0$, a contradiction. From now on we assume that
$g\neq\id$. Denote $\bar C:=\bar \Delta_{i_0}$ and let
$\delta_{i_0}=1-1/m$. Since $\delta(\bar S,\Dif{\bar S})=1$, there
are no divisors $E\neq \bar C$ of $\KKK(\bar S)$ with
$a(E,\Dif{\bar S})\le -6/7$. This gives us that $\bar S$ is smooth
at $\Sing{\bar C}$ whenever $\Sing{\bar C}\neq\emptyset$ (see
\cite[Lemma~9.1.8]{Pr-lect}). By our assumptions, $\bar Y$ is
singular along $\bar C$. Moreover, at the general point of $\bar
C$ we have an analytic isomorphism (see \cite[16.6]{Ut}):
\begin{equation}
\label{eq-2}
(\bar Y,\bar S,\bar C)\simeq
(\CC^3,\{x_3=0\},\{x_1-\text{axis}\})/\ZZ_m(0,1,q), \quad (m,q)=1.
\end{equation}

\begin{lemma}
\label{iso-gen}
Notation as above. Assume that $p_a(\bar C)\ge 1$. Then the map
$\bar Y\dashrightarrow Y$ is an isomorphism at the general point
of $\bar C$. Moreover, if $\bar P\in \bar C$ is a singular point,
then $\bar Y\dashrightarrow Y$ is an isomorphism at $\bar P$. In
particular, the proper transform $C$ of $\bar C$ is a curve with
$p_a(C)\ge 1$.
\end{lemma}

\subsection{}
\label{implication}
First we show that Lemma~\ref{iso-gen} implies Theorem~\ref{main}.
Assume $p_a(\bar C )\ge 1$. Clearly, $C\subset S$. By
Lemma~\ref{iso-gen}, $p_a(C)\ge 1$ and
\begin{equation}\label{number-of-singular-points}
\#\Sing{C}\ge\#\Sing{\bar C}
\end{equation}
From \eqref{eq-2}, we have $\Dif{S}=(1-1/m)C+(\text{other
terms})$.

\subsection{}
Consider the case when $g(S)$ is a point. By Lemma~\ref{first}, as
in the proof of Proposition~\ref{prop-main}, one can show that
there exists $1$, $2$, $3$, $4$, or $6$-complement of the form
$K_S+C+(\text{other terms})$ on $S$. By Adjunction, $p_a(C)=1$.
Therefore, $p_a(\bar C)=1$ and we have equality in
\eqref{number-of-singular-points}. Thus $S$ is smooth at
$\Sing{C}$ (whenever $\Sing{C}\neq\emptyset$). We have a
contradiction by Proposition~\ref{prop-main}.

\subsection{}
Consider the case when $g(S)$ is a curve. Note that $S$ is
rational (because $(\bar S, \Dif{\bar S})$ is a klt log del Pezzo,
see e.~g. \cite[Sect. 5.5]{Pr-lect}) and so is $g(S)$. Consider
the restriction $g_S\colon S\to g(S)$. Since $p_a(C)\ge 1$, $C$ is
not a section of $g_S$. Let $\ell$ be the general fiber of $g_S$.
Then $\ell\simeq\PP^1$ and
\[
2=-K_S\cdot \ell>\Dif{S}\cdot \ell\ge (1-1/m)C\cdot \ell\ge
\frac67C\cdot\ell.
\]
Thus $C\cdot \ell =2$ and $C$ is a $2$-section of $g_S$. Moreover,
\[
(\Dif{S}-(1-1/m)C)\cdot \ell<2-2(1-1/m)=2/m<1/2.
\]
Hence $\Dif{S}$ has no horizontal components other than $C$. Let
$P:=g(\ell)$, let $X'$ be a germ of a general hyperplane section
through $P$ and let $Y'\colon =g^{-1}(X')$. Consider the induced
(birational) contraction $g'\colon Y'\to X'$. Since singularities
of $X$ are isolated, $X'$ is smooth. By Bertini Theorem,
$K_{X'}+\ell$ is plt. Further, $Y'$ has exactly two singular
points $C\cap\ell$ and these points are analytically isomorphic to
$\CC^2/\ZZ_m(1,q)$ (see \eqref{eq-2}). This contradicts the
following lemma.

\begin{lemma}[see {\cite[Sect. 6]{Pr-lect}}]
Let $\phi\colon Y'\to X'\ni o'$ be a birational contraction of
surfaces and let $\ell:=\phi^{-1}(o')_{\red}$. Assume that
$K_{Y'}+\ell$ is plt and $X'\ni o'$ is smooth. Then $Y'$ has on
$\ell$ at most two singular points. Moreover, if $Y'$ has on
$\ell$ exactly two singular points, then these are of types
$\frac1{m_1}(1,q_1)$ and $\frac1{m_2}(1,q_2)$, where
$\gcd(m_i,q_i)=1$ and $\gcd(m_1,m_2)=1$.
\end{lemma}
\begin{proof}
We need only the second part of the lemma. So we omit the proof of
the first part. Let $\Sing{Y'}=\{P_1,P_2\}$. We use topological
arguments. Regard $Y'$ as a analytic germ along $\ell$. Since
$(X'\ni o')$ is smooth, $\pi_1(Y'\setminus\ell)\simeq
\pi_1(X'\setminus \{o'\})\simeq\pi_1(X')=\{1\}$. On the other
hand, for a sufficiently small neighborhood $Y'\supset U_i\ni P_i$
the map $\pi_1(U_i\cap \ell\setminus P_i)\to \pi_1(U_i\setminus
P_i)$ is surjective (see \cite[Proof of Th. 9.6]{K}). Using Van
Kampen's Theorem as in \cite[0.4.13.3]{Mo} one can show that
\[
\pi_1(Y'\setminus\{P_1,P_2\}) \simeq \langle \tau_1,
\tau_2\rangle/\{\tau_1^{m_1}=\tau_2^{m_2}=1,\ \tau_1\tau_2=1\}.
\]
This group is nontrivial if $\gcd(m_1,m_2)\neq 1$, a
contradiction.
\end{proof}

\begin{proof}[Proof of Lemma~{\ref{iso-gen}}]
The map $Y\dashrightarrow\bar Y$ is a composition of log flips:
\begin{equation}\label{eq-flips}
\begin{array}{lcccccccccc}
 Y=Y_0 &&&& Y_1 &&&& Y_2 &\cdots & Y_N=\bar Y \\
\downarrow\lefteqn{\scriptstyle{g}}& \searrow && \swarrow &&
\searrow &&\swarrow &&& \\
 X&& W_1 &&&& W_2 &&&&
\end{array}
\end{equation}
where every contraction $\searrow$ is $(K+S+(1+\ep)A)$-negative
and every $\swarrow$ is $(K+S+(1-\ep)A)$-negative.
Kawamata-Viehweg Vanishing \cite[Th. 1-2-5]{KMM} implies that
exceptional loci of these contractions are trees of smooth
rational curves \cite[Cor. 1.3]{Mo}. Thus Lemma~\ref{iso-gen} is
obvious if the curve $\bar C$ is nonrational. From now on we
assume that $\bar C$ is a (singular) rational curve.

\begin{lemma}
\label{contained}
Notation as above. Let $S_i$ be the proper transform of $S$ on
$Y_i$. If $f$ is not a flipping contraction, then $-S_i$ is ample
over $W_{i}$ for $i=1,\dots,N$. In particular, all nontrivial
fibers of $Y_i\to W_{i}$ are contained in $S_i$.
\end{lemma}
\begin{proof}
We claim that $-S_i$ is not nef over $Z$ for $i=1,\dots,N-1$.
Indeed, assume $\Exc(f)\neq f^{-1}(o)$. Take $o'\in f(\Exc(f))$,
$o'\neq o$ and let $\ell\subset g^{-1}(f^{-1}(o'))$ be any compact
irreducible curve. Clearly, $Y\dashrightarrow Y_i$ is an
isomorphism along $\ell$. Let $\ell_i$ be the proper transform of
$\ell$ on $Y_i$. Since $S\cdot \ell=0$, we have $S_i\cdot
\ell_i=0$. The curve $\ell_i$ cannot generate an extremal ray
(because extremal contractions on $Y_1,\dots,Y_{N-1}$ are
flipping). If $-S_i$ is nef over $Z$, then taking into account
that $\rho(Y_i/Z)=2$ we obtain $S_i\equiv 0$, a contradiction.
Thus we may assume that $\Exc(f)= f^{-1}(o)$ is a (prime) divisor.
Then the exceptional locus of $Y_i\to Z$ is compact. If $-S_i$ is
nef, this implies $S_i=\Exc(Y_i\to Z)$. Again we have a
contradiction because the proper transform of $\Exc(f)$ does not
coincide with $S_i$.

We prove our lemma by induction on $i$. It is easy to see that
$-S$ is ample over $W_0:=X$. The Mori cone $\NE(Y_i)$ is generated
by two extremal rays. Denote them by $R_i$ and $Q_i$, where $R_i$
(resp. $Q_i$) corresponds to the contraction $Y_i\to W_{i}$ (resp.
$Y_i\to W_{i+1}$). Suppose that our assertion holds on $Y_{i-1}$,
i.~e. $S_{i-1}\cdot R_{i-1}<0$. By our claim above, $S_{i-1}\cdot
Q_{i-1}> 0$ and after the flip $Y_{i-1}\dashrightarrow Y_{i}$ we
have $S_{i}\cdot R_{i}< 0$. This completes the proof of the lemma.
\end{proof}

\subsection{}
\label{proof}
Let us consider the case when $f$ is not a flipping contraction.
Let $C^{(i)}$ be the proper transform of $\bar C$ on $Y_i$. If
$p_a(C^{(i)})\ge 1$, then $Y_i\to W_i$ cannot contract $C^{(i)}$.
Thus $C^{(i+1)}$ is well defined. Now we need to show only that on
each step of \eqref{eq-flips} no flipping curves $\Exc(Y_i\to
W_i)$ pass through singular points of $C^{(i)}$. (Then
$Y_i\dashrightarrow Y_{i+1}$ is an isomorphism near singular
points of $C^{(i)}$ and we are done). By Lemma~\ref{contained} all
flipping curves $\Exc(Y_i\to W_i)$ are contained in $S_i$.
Therefore we can reduce problem in dimension two. The last claim
$\Exc(Y_i\to W_i)\cap \Sing{C^{(i)}}=\emptyset$ easily follows by
the lemma below.

\begin{lemma}
Let $\varphi\colon S\to \widehat{S}\ni \hat o$ be a birational
contraction of surfaces and let $\Delta=\sum \delta_i\Delta_i$ be
a boundary on $S$ such that $K_S+\Delta$ is klt and
$-(K_S+\Delta)$ is $\varphi$-ample. Put $\Theta:=\sum_{\delta_i\ge
6/7}\Delta_i$.  Assume that $\varphi$ does not contract components
of $\Theta$. Then $\Theta$ is smooth on $\varphi^{-1}(\hat
o)\setminus \Sing{S}$.
\end{lemma}
\begin{proof}
Assume the converse and let $P\in \Sing{\Theta}\cap
\bigl(\varphi^{-1}(\hat o)\setminus \Sing{S}\bigr)$. Let $\Gamma$
be a component of $\varphi^{-1}(\hat o)$ passing through $P$. Then
$\Gamma\simeq\PP^1$. There is an $n$-complement $\Delta^+=\sum
\delta_i^+\Delta_i$ of $K_S+\Delta$ near $\varphi^{-1}(\hat o)$
for $n\in\{1,2,3,4,6\}$ (see \cite[Th. 5.6]{Sh}, \cite[Cor.
19.10]{Ut}, or \cite[Sect. 6]{Pr-lect}). By the definition of
complements, $\delta_i^+\ge \min\{1, \down{(n+1)\delta_i}/n\}$ for
all $i$. In particular, $\delta_i^+=1$ whenever $\delta_i\ge 6/7$,
i.~e. $\Theta\le \Delta^+$. This means that $K_S+\Theta$ is lc.
Since $P\in \Sing{\Theta}$, $K_S+\Theta$ is not plt at $P$.
Therefore $\Theta=\Delta^+$ near $P$ and $\Gamma$ is not a
component of $\Delta^+$. We claim that $K_S+\Gamma$ is lc. Indeed,
$\Gamma$ is lc at $P$ (because both $\Gamma$ and $S$ are smooth at
$P$). Assume that $K_S+\Gamma$ is not lc at $Q\ne P$. Then
$K_S+(1-\varepsilon)\Gamma+\Delta^+$ is not lc at $P$ and $Q$ for
$0<\ep\ll 1$. This contradicts Connectedness Lemma \cite[5.7]{Sh},
\cite[17.4]{Ut}. Thus $K_S+\Gamma$ is lc and we can apply
Adjunction:
\[
(K_S+\Delta^++\Gamma)|_\Gamma\ge
(K_S+\Theta+\Gamma)|_\Gamma=K_\Gamma+\Diff{\Gamma}{\Theta}.
\]
Since $K_S+\Delta^+\equiv 0$ over $\widehat{S}$ and
$\Gamma\simeq\PP^1$, we have $\deg \Diff{\Gamma}{\Theta}< 2$. On
the other hand, the coefficient of $\Diff{\Gamma}{\Theta}$ at $P$
is equal to $(\Theta\cdot \Gamma)_P\ge 2$, a contradiction.
\end{proof}

\subsection{}
Finally let us consider the case when $f$ is a flipping
contraction. If $-S_i$ is ample over $W_{i}$ for $i=1,\dots,N$,
then we can use arguments of \ref{proof}. From now on we assume
that $S_I$ is nef over $W_I$ for some $1\le I\le N$.

Let $L$ be an effective divisor on $Z$ passing through $o$. Take
$c\in\QQ$ so that $K_X+cf^*L$ is lc but not klt. By Base Point
Free Theorem \cite[3-1-1]{KMM}, there is a member $M\in
\left|-n(K_X+cf^*L)\right|$ for some $n\in\NN$ such that
$K_X+cf^*L+\frac1n M$ is lc (but not klt). Thus, we may assume
$F=cf^*L+\frac1n M$. Let $K_Y+S+B'=g^*(K_X+cf^*L)$ be the crepant
pull-back. Write $B=B'+B''$, where $B',\ B''\ge 0$. Then
$-(K_Y+S+B')$ is nef over $Z$ and trivial on fibers of $g$. Run
$(K_Y+S+B')$-MMP over $Z$. Since $K_Y+S+B'\equiv -B''\neq 0$, this
$\QQ$-divisor cannot be nef until $S$ is not contracted. Therefore
after a number of flips we get a divisorial contraction:
\begin{equation}\label{eq-flips1}
Y\dashrightarrow Y_1\dashrightarrow \cdots\dashrightarrow Y_N=\bar
Y \dashrightarrow \cdots Y_{N'}\longrightarrow X'.
\end{equation}
Since $\rho(Y_i/Z)=2$ the cone $\NE(Y_i/Z)$ has exactly two
extremal rays. Hence the sequence \eqref{eq-flips} is contained in
\eqref{eq-flips1}.

\begin{claim}
$S_j$ is nef over $W_j$ and $-S_j$ is ample over $W_{j+1}$ for
$I\le j\le N'$.
\end{claim}
\begin{proof}
Clearly, $-S_I$ is ample over $W_{I+1}$ (because $S_I$ cannot be
nef over $Z$). After the flip $Y_I\dashrightarrow Y_{I+1}$ we have
that $S_{I+1}$ is ample over $W_{I+1}$. Continuing the process we
get our claim.
\end{proof}
Further, $X'$ has only terminal singularities. Indeed, $X'$ is
$\QQ$-factorial, $\rho(X'/Z)=1$ and $X'\to Z$ is an isomorphism in
codimension one. Therefore, one of the following holds:
\begin{enumerate}
\item
$-K_{X'}$ is ample over $Z$, then $X'\simeq X$;
\item
$K_{X'}$ is numerically trivial over $Z$, then so is $K_X$, a
contradiction;
\item
$K_{X'}$ is ample over $Z$, then $X\dashrightarrow X'$ is a flip
and $X'$ has only terminal singularities \cite[5-1-11]{KMM}.
\end{enumerate}
This shows also that $Y_{N'}\to X'$ is a plt blowup. Then we can
replace $X$ with $X'$ and apply arguments of \ref{proof}. This
finishes the proof of Theorem~\ref{main}.
\end{proof}

\subsubsection*{Concluding Remark} Shokurov's classification of
exceptional log del Pezzos with $\delta\ge 1$ uses reduction to
the case $\rho=1$. More precisely, this method uses the following
modifications: $\bar S\longleftarrow S^\bullet\longrightarrow
S^{\circ}$, where $S^\bullet\to \bar S$ is the blow up of all
divisors with discrepancy $a(E,\Dif{\bar S}\ge 6/7$, $S^\bullet\to
S^\circ$ is a sequence of some extremal contractions and
$\rho(S^\circ)=1$. Then all divisors with discrepancy $\ge 6/7$
are nonexceptional on $S^\circ$. In our case, a smooth elliptic
curve with coefficient $\ge 6/7$ on $S^\circ$ cannot be contracted
to a point on $\bar S$ (because the singularities of $\bar S$ are
rational). By Theorem~\ref{main} this case does not occur. The
situation in the case of a singular rational curve with
coefficient $\ge 6/7$ on $S^\circ$ which is contracted to a point
on $\bar S$ is more complicated. This case will be discussed
elsewhere.

\subsubsection*{Acknowledgements}
I thank V.~V.~Shokurov for encouraging me to write up this note. I
was working on this problem at Tokyo Institute of Technology
during my stay 1999-2000. I am very grateful to the staff of the
institute and especially to Professor S.~Ishii for hospitality.

\end{document}